\scrollmode
\documentclass[12pt]{amsart}
\usepackage{verbatim}
\usepackage{eucal}
\usepackage{amssymb}
\usepackage{eepic}
\usepackage{graphicx}
\usepackage{psfrag}
\usepackage{mathrsfs}
\usepackage[all]{xy}
\swapnumbers
\theoremstyle{plain}
\newtheorem{theorem}[subsection]{Theorem}
\newtheorem{proposition}[subsection]{Proposition}
\newtheorem{lemma}[subsection]{Lemma}

\theoremstyle{definition}
\newtheorem{definition}[subsection]{Definition}
\newtheorem{example}[subsection]{Example}

\theoremstyle{remark}
\newtheorem{remark}[subsection]{Remark}
\newcount\timehh\newcount\timemm
\timehh=\time 
\divide\timehh by 60 \timemm=\time
\newcommand{\draftauthor}[1]{\author{#1
    {
      --- \protect \protect\sc\today\ ---
      \ifnum\timehh<10 0\fi\number\timehh\,:\,\ifnum\timemm<10 0\fi\number\timemm
      \protect \, \, \protect \bf DRAFT
    }
  }
}
\newcommand{\R}{{\mathbb R}}

\newcommand{\C}{{\mathbb C}}
\newcommand{\Z}{{\mathbb Z}}
\newcommand{\Q}{{\mathbb Q}}

\newcommand{\V}{{\mathscr V}}
\newcommand{\B}{{\mathscr B}}
\newcommand{\HHH}{{\mathscr H}}
\newcommand{\SSS}{{\mathscr S}}
\newcommand{\UUU}{{\mathscr U}}
\newcommand{\OK}{{\mathscr O_{K}}}

\newcommand{\hp}{{\mathfrak H}}
\newcommand{\Vor}{Vorono\v{\i}}
\newcommand{\Int}{\operatorname{Int}}

\newcommand{\Dim}{\operatorname{dim}}
\newcommand{\hull}{\operatorname{hull}}
\newcommand{\Stab}{\operatorname{Stab}}
\newcommand{\Nor}{\operatorname{Nor}}
\newcommand{\cones}{\mathbf{C}}
\newcommand{\zyklen}{\mathbf{Z}}
\newcommand{\vv}{{\mathbf {v}}}

\newcommand{\xx}{{\mathbf {x}}}
\newcommand{\uu}{{\mathbf {y}}}
\newcommand{\yy}{{\mathbf {y}}}
\newcommand{\Hecke}{\mathbf{D}}

\newcommand{\sets}[1]{[\![#1]\!]}

\DeclareMathOperator*{\support}{supp}
\DeclareMathOperator*{\Comm}{Comm}

\CompileMatrices

\begin{document}

\title[Hecke operators and self-adjoint homogeneous cones]{Hecke
operators and $\Q $-groups associated to self-adjoint homogeneous cones}

\newif \ifdraft

\def \makeauthor{
\ifdraft
	\draftauthor{Paul Gunnells and Mark McConnell}
\else
\author{Paul E. Gunnells}
\address{Department of Mathematics\\
Columbia University\\
New York, NY  10027}
\email{gunnells@math.columbia.edu}

\author{Mark McConnell}
\address{Oklahoma State University\\
Department of Mathematics\\
401 Math Sciences\\
Stillwater, OK  74078--1058}
\email{mmcconn@math.okstate.edu}
\fi
}

\draftfalse

\date{November 20, 1998} 

\thanks{The first author was partially supported by NSF grant DMS
96--27870 and a Columbia University Faculty Research Grant.  The
second author was partially supported by NSF grant DMS--9704535.}


%
%
\begin{abstract}
Let $G$ be a reductive algebraic group associated to a self-adjoint
homogeneous cone defined over $\Q $, and let $\Gamma \subset G $ be a
appropriate neat arithmetic subgroup.  We present two algorithms to
compute the action of the Hecke operators on $H^{i}(\Gamma ;\Z )$ for
all $i$.  This simultaneously generalizes the modular symbol algorithm
of Ash-Rudolph \cite{ash.rudolph} to a larger class of groups, and
provides techniques to compute the Hecke-module structure of
previously inaccessible cohomology groups.
\end{abstract}

\maketitle

%
%
\section{Introduction}\label{introduction}
\subsection{}
Let $G$ be a reductive algebraic group defined over $\Q $, and let
$\Gamma \subset G (\Q )$ be a neat arithmetic subgroup.  The group
cohomology $H^{*} (\Gamma ;\Z )$ plays an important role in
contemporary number theory, through its connection with automorphic
forms and representations of the absolute Galois group.
See~\cite{ash.galois} for an introduction to these ideas.  This
relationship is revealed in part through the action of the \emph{Hecke
operators} on the complex cohomology $H^{*} (\Gamma ;\C)$.  These are
endomorphisms induced from a family of correspondences associated to
the pair $(\Gamma ,G (\Q ))$ (\S\ref{hecke.operators}); the arithmetic
nature of the cohomology is contained in the eigenvalues of these
linear maps.

To compute the Hecke action in certain cases, one may use the
\emph{modular symbol algorithm.}  One begins with a finite cell
complex that computes $H^{*} (\Gamma ;\Z )$, and applies the Hecke
operators to an easily understood set of dual homology classes, the
\emph{modular symbols} (\S\ref{compare}).  There is a finite set of
modular symbols, distinguished by the choice of finite cell complex,
that spans the homology classes.  The modular symbol algorithm enables
one to write the Hecke-image of a modular symbol as a sum of symbols
taken from the finite spanning set.

Several research groups have used this technique to produce many
corroborative examples of the ``Langlands philosophy''
\cite{exp.ind, apt, crem, geemen.top}.  However, this
technique has some shortcomings:
\begin{itemize}
\item The group $G$ must be either a linear group ($SL_{n}$ or
$GL_{n}$)~\cite{ash.rudolph} or a symplectic group ($Sp_{2n}$)
\cite{gunn.symp}. 
\item The group $\Gamma $ must be associated to a euclidean domain.
For example, $\Gamma \subset SL_{n} (\Z )$ can be studied, but not
$\Gamma \subset SL_{n} (\OK )$, where $\OK $ is the ring of integers
in an algebraic number field $K/\Q $ with class number $>1$.
\item The modular symbol algorithm enables computation of the
Hecke action only on $H^{d}(\Gamma ;\Z )$, where $d$ is the
\emph{cohomological dimension} of $\Gamma $.  This is the smallest
integer $d$ such that $H^{*} (\Gamma ;M)=0$ for $*>d$ and any $\Z \Gamma
$-module $M$.
\end{itemize}

These limitations, particularly the last two, are real obstacles to
continuing the experimental work cited above.  For example:

\begin{itemize}
\item One is interested in $H^{2}$ of \emph{Bianchi groups},
especially the cuspidal classes~\cite{grune.sch}.  These arithmetic
groups have the form $\Gamma \subset SL_{2} (\OK )$, where $K/\Q $ is
an imaginary quadratric extension.  At present no systematic
``motivic'' explanation of the Hecke eigenvalues has been conjectured,
although in some cases one can make the connection with the absolute
Galois group~\cite{taylor}.
\item Work of Avner Ash and David Ginzburg suggests that one can find
new automorphic $L$-functions by integrating rational cuspidal
cohomology classes~$\alpha$ for $G = GL_4$ over $(2,1,1)$-modular
symbols, which are certain submanifolds of the associated locally
symmetric space.  Here $\alpha\in H^5$.  If the $(2,1,1)$-modular
symbols span a dual space to the cuspidal $H^5$ under this pairing,
then the $L$-functions would exist and be non-zero.  For $\Gamma_0(N)
\subset SL_4(\Z)$, Ash and the second author have computed
$H^5(\Gamma_0(N))$ for a range of prime levels~$N$.  Computing the
Hecke action on these groups will tell us, at least conjecturally,
which classes are cuspidal, and which are lifts from smaller groups
like $Sp_4$ or $O_4$.  The modular symbol algorithm cannot be used
because it works only in degree~6, not~5.
\end{itemize}

\subsection{}
In this paper we overcome these obstacles for a special class of
arithmetic groups, the groups $\Gamma \subset G (\Q )$ such that $G$
is the automorphism group of a \emph{self-adjoint homogeneous cone}
with a linear structure compatible with the $\Gamma $-action
(\S\ref{bk.cones}).  The real groups involved are $G (\R )= GL_{n} (\R
)$ and $GL_{n} (\C)$, as well as more exotic examples
(\S\ref{rationality}).  The relevant arithmetic groups include $\Gamma
\subset G(R)$, where $R$ is
\begin{itemize}
\item $\Z $,
\item the ring of integers in a totally real field, and 
\item the ring of integers in a $CM$ field.
\end{itemize}
We present two algorithms (in Theorems~\ref{thm1},~\ref{thm2}, and
in Theorem~\ref{thm3}) for the computation of the Hecke action.  

\subsection{}\label{compare}
In the remainder of this introduction, we compare our algorithms with
the modular symbol algorithm of Ash-Rudolph~\cite{ash.rudolph}.  For
concreteness, we work with the simplest example with more than one
interesting cohomology group, that of $SL_{3}$.  First we establish
notation.

Let $G$ be the split form of $SL_{3}$, so that $G (\Q )=SL_{3} (\Q )$.
Let $V$ be the $\R $-vector space of all $3\times 3$ real symmetric
matrices, and let $C\subset V$ be the open cone of positive-definite
matrices.  Then $G (\R ) = SL_{3} (\R )$ acts on $C$ by $c\mapsto
g\cdot c\cdot g^{t}$, and the stabilizer of $c$ is isomorphic to
$SO_{3} (\R )$, the maximal compact subgroup of $SL_{3}(\R )$.  The
group $\R ^{>0}$ acts on $C$ by homotheties, and we denote the
quotient by $X$.  We have an isomorphism $SL_{3} (\R )/SO_{3} (\R
)\xrightarrow{\sim} X$ given by $gK\mapsto \R ^{>0}gg^{t}$, and $X$ is
a smooth noncompact manifold of real dimension $5$.

Let $L\subset V$ be the lattice of integral symmetric matrices, and
let $\Gamma _{L}= SL_{3} (\Z )$ be the stabilizer of $L$.  We fix a
neat $\Gamma \subset \Gamma _{L}$.  The space $X$ is contractible, so
we may identify the group cohomology $H^{*} (\Gamma ;\Z )$ with $H^{*}
(\Gamma\backslash X;\Z )$.  The space $\Gamma \backslash X$ is
$5$-dimensional, and the cohomological dimension of $\Gamma $ is $3$.

\subsection{}\label{def.of.ms}
We may use modular symbols to study $H^{3} (\Gamma )$.  Let $\bar C$
be the closure of $C$.  Any nonzero primitive vector $v\in \Z ^{3}$
determines a rank-one semidefinite quadratic form $q (v)\in \bar
C\smallsetminus C$ as follows: if we write $v= (a,b,c)^{t}\in \Z ^{3}$,
then
\[
q (v) := 
\left(\begin{array}{c}
a\\
b\\
c
\end{array} \right) \left(\begin{array}{ccc}
a&b&c
\end{array} \right).
\]
The rays $\left\{\R ^{>0}q (v)\mid v\in \Z ^{3}\smallsetminus \{0 \}
\right\}$ are called the \emph{cusps} of $C$.  

Let $\vv =
(v_{1},v_{2},v_{3}) $ be an ordered triple of distinct nonzero
primitive integral vectors.  Then if $\vv$ is linearly independent, 
it determines an open oriented $3$-cone
$\sigma (\vv ) \subset C$, by
\[
\sigma (\vv) = \Bigl\{\sum \rho _{i}q(v_{i}) \Bigm |  \rho _{i}>0 \Bigr\}.
\]
Under the composition 
\[
C\longrightarrow X\longrightarrow \Gamma \backslash X,
\]
this cone is taken to an open submanifold of $\Gamma \backslash X$.
The closure of this submanifold in $\Gamma \backslash \bar X$, the
Borel-Serre compactification of $\Gamma \backslash X$
(\S\ref{compactifications}), determines a class $[\vv ]\in H_{2}
(\Gamma \backslash \bar X,\partial (\Gamma \backslash \bar X))$.  Such
a class is by definition a \emph{modular symbol}.  Via Lefschetz
duality we may identify $[\vv ]$ with a class in $H^{3}(\Gamma
\backslash \bar X) = H^{3} (\Gamma ;\Z )$.  It can be shown that the
duals of the modular symbols span $H^{3} (\Gamma ;\Z )$.  However,
there are infinitely many of them.

If $\det \vv =\pm 1$, then $[\vv ]$ is called a \emph{unimodular
symbol}.  There are only finitely many unimodular symbols modulo
$\Gamma $, and one can show that their duals span $H^{3}$ as
follows~\cite{ash.rudolph}.  Suppose $|\det \vv |>1$.  Using the
euclidean algorithm, one can construct a nonzero $w\in \Z ^{3}$ such
that
\[
0\leq |\det (w,v_{i},v_{j})| < |\det \vv |
\]
for any $1\leq i<j\leq 3$.  Since modular symbols satisfy the relation 
\[
[v_{1},v_{2},v_{3}] = [w,v_{2},v_{3}] - [w,v_{1},v_{3}] + [w,v_{1},v_{2}],
\]
and since $[\vv ]=0$ if $\det \vv = 0$, by iterating one can write a
modular symbol as a sum of unimodular symbols.  This is the
\emph{modular symbol algorithm}.

The image of a unimodular symbol $\alpha $ under a Hecke operator $T$
is a finite sum of modular symbols, which in general are not
unimodular.  Using the modular symbol algorithm, we may write $T
(\alpha )$ as a sum of unimodular symbols.  Thus we may compute the
Hecke action on $H^{3} (\Gamma ;\Z )$ using the finite spanning set
given by the unimodular symbols.

\subsection{}\label{whats.new}
Now we describe our approach to computing the Hecke action.  We begin
by shifting attention from the tuple $\vv $ to the cone $\sigma (\vv
)$.  Instead of using $\det \vv $ as a measure of non-unimodularity,
we use the relative position of $\sigma (\vv)$ with respect to a
distinguished collection of cones in $\bar C$, the set $\V$ of
\emph{\Vor \ cones} (\S\ref{gam.admis.section}).  In this example,
$\V$ is the $SL_{3} (\Z )$-orbit of the closed cone generated by the
six cusps $\R ^{>0}q (e_{i})$ and $\R ^{>0}q (e_{i}-e_{j})$, where
$\{e_{i} \}\subset \Z ^{3}$ is the standard basis, and $1\leq i<j\leq
3$.  The intersection of each of these cones with $C$ is a weak
fundamental domain for $SL_{3} (\Z )$, in the following sense: any
$x\in C$ meets a non-empty finite subset of $\V $.

One may use $\V$ for cohomology calculations as follows.  First we let
$\cones^{R}_{*}$ be the complex over $\Z $ generated by \emph{all}
simplicial rational cones in $\bar C$ of all dimensions, with the
obvious boundary maps (\S\ref{realization}).  This complex maps by
linear projection to the singular chain complex of a certain Satake
compactification $\Gamma \backslash \tilde{X} $ of $\Gamma \backslash
X$, and is surjective on homology.  The subcomplex $\cones^{V} _{*}$
generated by the \Vor \ cones is a finite complex mod $\Gamma $, and
also maps surjectively to homology.  The relative homology $H_{*}
(\Gamma \backslash \tilde{X}, \Gamma \backslash \partial \tilde{X})$
can then be identified with the cohomology $H^{5-*} (\Gamma \backslash
X)= H^{5-*} (\Gamma )$ (Proposition~\ref{comparing.coho}).

The complex $\cones^{R} _{*}$ is a variant of a complex that first
appeared in the literature in~\cite{lee.szczarba} and
\cite{ash.unstable}\footnote{This latter complex is called the
\emph{sharbly complex}, in honor of the authors of
\cite{lee.szczarba}.  The name is due to Lee Rudolph.}, and it plays
the role of the modular symbols---it is infinite mod $\Gamma $ and is
preserved by the Hecke operators.  The \Vor \ subcomplex $\cones^{V}
_{*}$, on the other hand, plays the role of the unimodular
symbols---it is finite mod $\Gamma $, yet is not preserved by the
Hecke operators.  Thus to compute the Hecke action, we must show how
to write any cycle in $\cones^{R}_{*}$ that is a Hecke image of a \Vor
\ cycle as a sum of \Vor \ cycles.  Hence we have replaced the
algebraic problem of reducing a determinant with the geometric problem
of moving a cycle built of ``generic'' cones into a cycle supported on
\Vor \ cones.

\subsection{}
\begin{figure}[ht]
\begin{center}
\psfrag{sf}{$F (\vv )$}
\psfrag{s}{$\sigma (\vv )$}
\psfrag{v1}{$v_{1}$}
\psfrag{v2}{$v_{2}$}
\psfrag{v3}{$v_{3}$}
\psfrag{w1}{$w_{1}$}
\psfrag{w2}{$w_{2}$}
\psfrag{w3}{$w_{3}$}
\psfrag{w4}{$w_{4}$}
\includegraphics{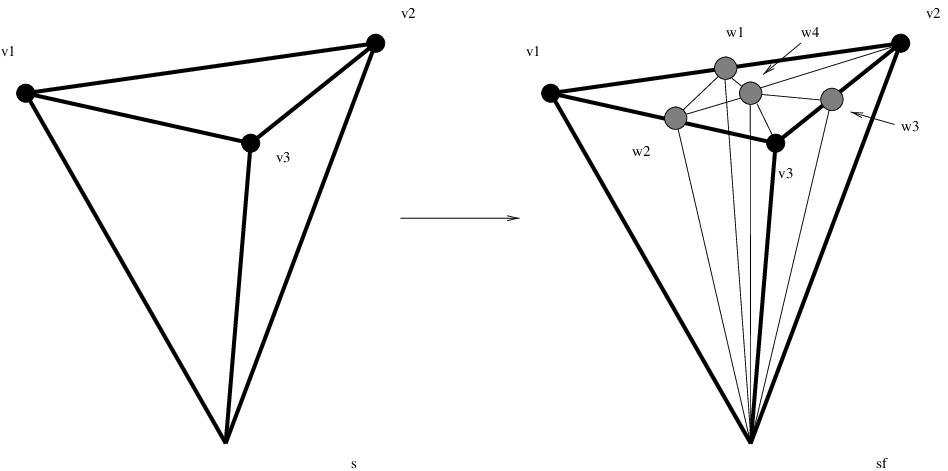}
\end{center}
\caption{\label{canonicalfig}}
\end{figure}

We give two techniques to do this.  Our first technique
(Theorem~\ref{thm1}) replaces $\sigma (\vv )$ with a set of cones $F
(\vv )$ constructed by refining the intersection of $\sigma (\vv )$
with the \Vor \ cones.  Then for each $1$-dimensional cone of $F (\vv
)$ we choose a cusp $w$.  These cusps can be combined with the
original $v\in \vv $ into tuples, and this yields a relation in
homology.  For example, in Figure~\ref{canonicalfig}, we show a cone
spanned by $\vv = (v_{1}, v_{2}, v_{3})$, which represents a class in
$H^{3} (\Gamma \backslash X)$.  This class is equal to the sum (with
appropriate orientations) of the classes corresponding to $(v_{1},
w_{1}, w_{2})$, $(w_{1}, w_{2}, w_{4})$, \dots, $(v_{3},w_{3},w_{4})$.
The new cones are \Vor \ cones, so we have reached our goal: up to
homology, we may replace $\sigma (\vv )$ with a sum of \Vor \ cycles.

In general, a cycle $\xi \in \cones ^{V}_{*}$ will be a sum of cones
with vanishing boundary mod $\Gamma $, and its image under a Hecke
operator $T (\xi )$ will be a similar object in $\cones ^{R}_{*}$.  We
show that the new cusps can be chosen $\Gamma $-equivariantly over $T
(\xi) $, which ensures that the result is a cycle
(Theorem~\ref{thm2}).  We further show that the cusps can be chosen so
that the resulting cones are \Vor \ cones.

Our first technique has the disadvantage that the fans $F (\vv )$ are
difficult to compute in practice.  In our second technique, we
circumvent this problem by using the \emph{\Vor \ reduction algorithm}
and a \emph{sufficiently fine decomposition} of $\sigma (\vv )$
(\S\ref{second}).  The former is an algorithm that computes which \Vor
\ cone contains a given point of $C$, and the latter is
decomposition of $\sigma (\vv )$ into cones that are small enough to
construct the homology (Theorem~\ref{thm3}).  We show how a
sufficiently fine decomposition may be used to transform a generic
cycle into a cycle built of \Vor \ cones (Theorem~\ref{thm4}).

\subsection{Acknowledgements} 
We thank Avner Ash and Bob MacPherson for encouragement and support.
The first author thanks Oklahoma State University for its hospitality.
We are both grateful to the NSF for support.

%
%
\section{Background}\label{background}
In this section we recall facts about self-adjoint homogeneous cones
and Hecke correspondences.  For more details, the reader may consult
\cite{ash, cones} for cones and~\cite{krieg, shimura} for
Hecke correspondences.  \S\S\ref{bk.cones}--\ref{arithmetic.gps}
closely follow Ash~\cite{ash}.

\subsection{}\label{bk.cones}
Let $V$ be an $\R $-vector space defined over $\Q $, and
let $C\subset V$ be an open cone.  That is, $C$ contains no straight
line, and $C$ is closed under homotheties: if $x\in C$ and $\lambda
\in \R^{>0}$, then $\lambda x\in C$.  The cone $C$ is
called {\em self-adjoint} if there exists a scalar product $\langle\,
,\, \rangle$ on $C$ such that
$$C = \left\{\,x\in V \,\bigm |\, \langle y,x\rangle > 0 \quad \hbox{for $y\in \bar C\smallsetminus \{0 \}$} \,\right\}. $$
Such a cone is necessarily convex, as is $\bar C$.

Let $G$ denote the connected component of the identity of the linear
automorphism group of $C$, i.e. $G=\left\{\,g\in GL(V)\mid gC = C\,
\right\}^{0}$.  The cone $C$ is called {\em homogeneous} if $G$ acts
transitively on $C$.  If $K$ denotes the isotropy group of a given 
point in $C$, then we may identify $C$ with $G/K$.  The
self-adjointness of $C$ implies that $G$ is reductive and that $C$
modulo homotheties is a riemannian symmetric space.  We denote this
symmetric space throughout by $X$, and always write $N=\Dim X$.

\subsection{}\label{rationality}
We assume that all the above notions are compatible with the $\Q
$-structure on $V$.  That is, as a subgroup of $GL(V)$, $G$ is defined
by rational equations, and the scalar product 
is defined over $\Q $.  These rationality conditions place the
following restrictions on $V$.

Recall that a \emph{Jordan algebra} is a finite-dimensional
algebra $J$ satisfying 
\begin{enumerate}
\item $ab=ba$,
\item $a^{2} (ba) = (a^{2}b)a$,
\end{enumerate}
for all $a,b\in J$.  In general $J$ is not associative.  If $J$ is
defined over $\R $, we say that $J$ is \emph{euclidean} if $a^{2}
+ b^{2} = 0$ implies $a=b=0$. 

Fix a basepoint $p\in C (\Q):= C\cap V (\Q )$.  Then the rationality
conditions on $G$ hold if and only if $V$ can be given the structure
of a euclidean Jordan algebra defined over $\Q $ with identity $p$
such that $C$ is the set of invertible squares in $V$.  This implies
that the group of real points $G (\R )$ must be isomorphic to a
product of the following groups \cite[p. 97]{cones}:
\begin{enumerate}
\item $GL_{n} (\R )$.
\item $GL_{n} (\C )$.
\item $GL_{n} (\mathbb H )$.
\item $O(1,n-1)\times \R ^{\times }$.
\item The noncompact Lie group with Lie algebra $\mathfrak e_{6
(-26)}\oplus \R $ (This is the group of
collineations of the projective plane over the Cayley numbers
\cite[p. 46]{tits}.) 
\end{enumerate}
In each case $V$ is a set of hermitian symmetric matricies.  In other
words, $V$ is the set of $n\times n$ matrices $A$ over an appropriate $\R
$-algebra with involution $\tau $, such that
$A^{t}= A^{\tau }$.  The cone $C$ is then the subset of
``positive-definite'' matrices in an appropriate sense.  For details
we refer to \cite[Ch. V]{cones}.

\subsection{}\label{const.bdy.comps}
Let $H$ be a hyperplane in $V$.  We say that $H$ is a {\em supporting
hyperplane} of $C$ if $H$ is rational and $H\cap C=\varnothing $ but
$H\cap \bar C \not =\varnothing $.  Since $\bar C$ is convex, these
conditions imply that $\bar C$ lies entirely in one of the two closed
half-spaces determined by $H$.

Given a rational supporting hyperplane $H$ of $C$, let $C' =
\Int(H\cap \bar C)$.  (Throughout, $\Int(A )$ is the interior of $A$
in its linear span.)  Then $C'$ is called a {\em rational boundary
component}, and is a self-adjoint homogeneous cone of smaller
dimension than $C$.  By convention, we also say that $C$ is a
(improper) rational boundary component.  Let $\tilde{C}$ be the union
of $C$ and all its proper rational boundary components.

Here is characterization of the rational boundary components
in terms of the Jordan algebra structure on $V$.  If $e\in V (\Q )$
satisfies $e^{2}=e$, then we call $e$ a \emph{rational idempotent},
and write $V (e)=\{x\in V \mid xe = x \}$.  The subspace $V (e)$ is a
euclidean sub Jordan algebra defined over $\Q $.  Then any
rational boundary component arises as the subset of invertible squares
$C (e)\subset V (e)$ for some choice of $e$.  

Two idempotents $e$ and $f$ are called \emph{orthogonal} if $ef=0$.  A
rational idempotent is called \emph{minimal} if it cannot be written
as the sum of two orthogonal nonzero idempotents.  Any rational
idempotent $e$ can be written as a sum of mutually orthogonal minimal
idempotents, and the number needed is an invariant of $e$.  This
number is called the \emph{$\Q $-rank} of $e$.  By definition, the
$\Q$-rank of a rational boundary component $C (e)$ is the $\Q $-rank
of $e$.  For any nonnegative integer $k$, let $C (k)$ denote the union
of all rational boundary components of $\Q $-rank $\leq k$.

\begin{definition}\label{cusps}
The {\em cusps} of $C$ are the rank-one rational boundary
components of $C$.  The set of cusps is denoted $\Xi (C)$.
\end{definition}

We will always use $n$ to denote the $\Q $-rank of $C$ itself.  Note
that $C (1)=\Xi (C)$, and that $C (k)$ is disjoint from $C$ if $k<n$.
If $k\geq n$, then $C (k)=\tilde{C}$.

\begin{remark}
Because of the homotheties, the $\Q $-rank of $C$ is \emph{not} the
same as the $\Q $-rank of the algebraic group $G$.  If $n$ is the $\Q
$-rank of $C$, then $n-1$ is the $\Q $-rank of $G$.  For example, when 
$G=SL_{n}$, the $\Q $-rank of $C$ is $n$.
\end{remark}

\subsection{}\label{arithmetic.gps}
Let $L\subset V(\Q )$ be a lattice, i.e. a discrete subgroup of $V(\Q
)$ such that $L\otimes \Q =V(\Q )$.  Let $\Gamma_{L}$ denote the
subgroup of $G(\Q )$ carrying $L$ onto itself.  An {\em arithmetic
subgroup } of $G$ is a discrete subgroup commensurable with
$\Gamma_{L}$ for some $L$.  Any neat arithmetic group $\Gamma \subset
\Gamma_{L}$ of finite index acts properly discontinuously and freely
on $C$.  Thus the quotient $\Gamma \backslash C$ is an
Eilenberg-Mac~Lane space for $\Gamma $, and the group cohomology
$H^{*}(\Gamma )$ is $ H^{*}(\Gamma\backslash C )$.  In fact, since
homotheties commute with the action of $\Gamma $, we may pass to the
symmetric space $X$, and compute $H^{*}(\Gamma \backslash X)$ instead.

\subsection{}\label{hecke.operators}
Fix a neat arithmetic group $\Gamma $.  Given $g\in G (\Q )$, let
$\Gamma ^{g} = g^{-1}\Gamma g$ and $\Gamma ' = \Gamma \cap \Gamma^{g}$.  Let $\Comm(\Gamma)$
be the commensurator of $\Gamma $.  This is the subgroup of $G (\Q
)$ defined by
\[
\Comm(\Gamma) := \Bigl\{g\in G (\Q )\,\Bigm |\,[\Gamma : \Gamma'] \quad \hbox{and} \quad   [\Gamma^g: \Gamma'] < \infty \Bigr\}.
\]  

For any $g\in \Comm (\Gamma )$, the inclusions $\Gamma '\rightarrow
\Gamma $ and $\Gamma '\rightarrow \Gamma^ g$ determine a diagram
\[
\vbox{\xymatrix{&{\Gamma '\backslash X}\ar[dr]^{t}\ar[dl]_{s}&\\
          {\Gamma \backslash X}&&{\Gamma \backslash X}}}
\]
Here $s(\Gamma 'x) = \Gamma x$ and $t$ is the composition of $\Gamma '
x \mapsto \Gamma ^{g}x$ with left multiplication by $g$.  This diagram
is called the {\em Hecke correspondence} associated to $g$.  It can be
shown that, up to isomorphism, the Hecke correspondence depends only
on the double coset $\Gamma g\Gamma $.  Furthermore, it can also be
shown that Hecke correspondences extend naturally to the rational
boundary components, inducing a correspondence\footnote{Here one takes
the Satake topology for the quotient (cf. \S\ref{compactifications}).}
\[
\vbox{\xymatrix{
&{\Gamma '\backslash \tilde{X}}\ar[dr]^{\tilde{t}}\ar[dl]_{\tilde{s}}&\\
{\Gamma \backslash \tilde{X}}&&{\Gamma \backslash \tilde{X}}}}
\]
Because the maps $s$ and $t$ are proper, we obtain a map on
cohomology
\[
T_{g} := t_{*}s^{*}\colon H^{*}(\Gamma \backslash X;\Z
)\rightarrow H^{*}(\Gamma \backslash X;\Z). 
\]
 
This is called the {\em Hecke operator} associated to $g$.
We let ${\HHH}_{\Gamma }$ be the $\Z $-algebra
generated by the Hecke operators, with product given by composition.

%
%
\section{Rational Polyhedral Cycles}\label{rpc}
In this section we develop some tools to compute the cohomology of
$\Gamma $.  The main result is a correspondence between classes in
$H^{*}(\Gamma ;\Z )$ and homology classes in a chain complex built
from rational polyhedral cones.  Our construction relies on the
reduction theory for self-adjoint homogeneous cones
(\S\ref{gam.admis.section}) due to Ash~\cite{ash}.

\subsection{}\label{retract}
Let $A\subset V (\Q )$ be a finite set of nonzero points.  The closed
convex hull $\sigma $ of the rays $\left\{\R ^{\geq 0}x \bigm| x\in A
\right\}$ is called a {\em rational polyhedral cone}.  We say that
$\sigma $ is a $d$-cone if it has dimension $d$.  The rays through the
vertices of the convex hull of $A$ are called the {\em spanning rays}
of $\sigma $.  We denote the set of spanning rays by $R(\sigma )$.
The group $G (\Q )$ acts naturally on the set of rational polyhedral
cones, and we denote the action by a dot: $\sigma \mapsto g\cdot
\sigma $.

Let $\xx$ be a finite tuple of nonzero points in $V (\Q )$.  Then a
\emph{pointed rational polyhedral cone} $\sigma (\xx )$ is the data of
the cone generated by the $x\in \xx $ and the tuple $\xx $.  There is
a natural $G (\Q )$-action on pointed rational polyhedral cones given
by the action on tuples, and this action agrees with the action on the
underlying cones.  We will sometimes suppress the tuple $\xx $ from
the notation.

We say that $\sigma (\xx )$ is \emph{simplicial} if every subset of $R
(\sigma (\xx )) $ spans a face of $\sigma (\xx )$.  This means that
$\sigma (\xx )$ is equivalent to a $d$-simplex (except it is missing a
facet at infinity in $V$), where $d = \#R (\sigma (\xx ))$.

\subsection{}\label{fan.stuff} 
A collection of cones $F$ is called a \emph{fan} 
if it satisfies the following:
\begin{enumerate}
\item Any face $\tau $ of $\sigma \in F$ is also a member of $F$.
\item If $\sigma,\tau \in F$, then $\sigma
\cap \tau $ is a common face of $\sigma $ and $\tau $.
\end{enumerate}
The subsets $F_{d} := \left\{\sigma \in F_{d}\mid \dim \sigma \leq d
\right\}$ provide a filtration of $F$ by subfans.  We do not require
that $F$ be either finite or locally finite.

\subsection{}\label{gam.admis.section}
Now we want to partition $C$ into convex subsets using a
rational polyhedral fan, in a manner compatible with
the $\Gamma_{L}$-action.  This requires some care, as $C$ is open
and any rational polyhedral cone is closed.

Let $\Gamma \subset G$ be an arithmetic group.
A fan $F$ is called a
{\em $\Gamma $-admissible decomposition} of $C$ when the following
hold \cite[p. 72]{ash}:
\begin{enumerate}
\item Each $\sigma\in F$ is the span of a finite number of
rational vertices, and is a subset of $\bar C$.
\item For any $\sigma\in F$ and any $\gamma \in \Gamma $,
we have  $\gamma \cdot \sigma \in F$.
\item The set $\Gamma \backslash F$ is finite.
\item $C = \bigcup _{\sigma \in F}(\sigma\cap C)$.
\end{enumerate}
Note that a $\Gamma $-admissible decomposition descends mod
homotheties to a
decomposition of $X$ into open cells.

Here is a technique to construct $\Gamma $-admissible decompositions.
It originated with \Vor \cite{voronoi1}, and was generalized by Ash to
all self-adjoint homogeneous cones \cite[Ch. II]{amrt}.  Let $L'$ be
$L\smallsetminus \{0 \}$.

\begin{definition}\label{vor.poly}
The {\em \Vor \ polyhedron} $\Pi $ is the closed convex hull of
$L'\cap \,\Xi (C)$.
\end{definition}

According to the Lemma of \cite[p. 75]{ash} and its proof, the
polyhedron $\Pi $ has vertices only in $\Xi (C)$, and each face of
$\Pi $ is the convex hull of a finite set of points.  Clearly $\Pi
\subset \bar C$, since $\Xi (C)\subset \bar C$.

\begin{theorem}
\cite[p. 143]{amrt} The cones over the faces of $\Pi $ form a
$\Gamma $-admissible decomposition of $C$.
\end{theorem}

\begin{proposition}\label{gam-eq-for-Pi}
\cite[p. 75]{ash} Suppose that $\Gamma $ is neat.  Then the fan $F$
of cones over the faces of $\Pi $ may be refined, without adding new
one-dimensional fans, to a $\Gamma $-invariant simplicial fan $F'$.
\end{proposition}

\begin{proof}
We present the proof, since it is the prototype for similar arguments
later.  Every cone can be subdivided into simplicial cones without
adding new rays; we must prove that this can be done over all of $F$
in a $\Gamma $-invariant way.  We construct refinements $F'_{k}$ of
$F_{k}$ by induction on the dimension.  Our $F'$ will be $F'_{N}$,
where $N$ is the dimension of $V$.

To begin, any 1-cone is simplicial.  Hence we set $F'_{1}= F_{1}$.

Now suppose that $F_{k}$ has been $\Gamma $-equivariantly refined to
$F_{k}'$.  The set of cones $F_{k+1}\smallsetminus F_{k}$ is finite
modulo $\Gamma $, and we may choose a set $T$ of representatives of
these orbits.  Let $\sigma \in T$.  Let $\Nor (\sigma )$ be the
subgroup of $\Gamma $ preserving $\sigma $ as a set, and let $\Stab
(\sigma )$ be the subgroup of $\Gamma $ preserving each spanning ray
of $\sigma $.  Elements of $\Stab (\sigma )$ actually fix $\sigma $
pointwise, since they must fix the primitive generator of each
spanning ray in $L'$.  The quotient group $\Nor (\sigma )/\Stab
(\sigma )$ is finite, since it is a permutation group on the primitive
generators of the spanning rays of $\sigma $.  Because $\Gamma $ is
neat, this finite group is trivial, and hence $\Nor (\sigma ) = \Stab
(\sigma ) $.

Now choose any subdivision of $\sigma $ into simplicial cones without
adding new rays.  Since $\Stab (\sigma )$ preserves this subdivision,
so does $\Nor (\sigma )$.  Hence we may use the $\Gamma $-action to
carry the subdivisions for the $\sigma \in T$ to all of $F_{k+1}$ in a
well-defined way.
\end{proof}

For the remainder of this paper, we assume $\Gamma $ is neat, and we
fix a $\Gamma $-invariant simplicial refinement of the cones over the
faces of $\Pi $.  We call the resulting simplicial cones the {\em \Vor
\ cones}, and denote the fan of \Vor \ cones by $\V $.

We can make any \Vor \ cone $\sigma $ into a pointed cone
$\sigma (\xx )$ by choosing an ordering of $R (\sigma )$.  The
rational points comprising $\xx $ are the corresponding vertices of
$\Pi $.

\subsection{} 
As in \S\ref{const.bdy.comps}, any rational boundary component $C'$ of
$C$ is itself a self-adjoint homogeneous cone with $\Q $-rank less
than that of $C$.  Let $\Xi (C') = \Xi (C)\cap C'$, the set of cusps
in $C'$.  The \Vor \ fan $\V$ is compatible with the rational
boundary components in the following sense.

\begin{proposition}\label{bdy}
Let $C'$ be a rational boundary component of $C$ and let $\Pi '$ be
$C' \cap \Pi$.  Then $\Pi '$ is a face of $\Pi $, and is the
convex hull of $L'\cap \Xi ( C')$.
\end{proposition}

\begin{proof}
Let $H$ be a supporting hyperplane of $\bar C$ cutting out $C'$.
Section~\ref{const.bdy.comps} implies $H$ meets $\Xi (C')$, so $H\cap
\Pi \not = \varnothing$.  Since $\Pi \subseteq \bar C$, $H$ is a
supporting hyperplane of $\Pi $.  We have $\Pi ' = C'\cap \Pi = H\cap
\bar C \cap \Pi = H\cap \Pi $.  Thus $\Pi '$ is cut out from $\Pi $ by
a supporting hyperplane; by definition, this means $\Pi '$ is a face
of $\Pi $.

Let $\hull (S)$ denote the convex hull of a set $S$.  We have $\Pi ' =
H \cap \Pi = H \cap \hull(L'\cap \Xi (C))$.  Since $\Xi (C)$ lies
entirely in one of the closed half-spaces determined by $H$, the
latter set equals $\hull (H\cap L'\cap \Xi (C))$.  In turn, this is
$\hull (L'\cap H\cap \bar C \cap \Xi
(C)) = \hull (L' \cap C'\cap \Xi (C)) = \hull (L'\cap \Xi (C'))$.
\end{proof}

\subsection{}\label{compactifications}
Let $\bar X$ be the bordification of the symmetric space $X$
constructed by Borel-Serre~\cite{borel.serre}.  Let $Y=\Gamma
\backslash X$ and $\bar Y=\Gamma \backslash \bar X$.  Let $\partial
\bar Y = \bar Y \smallsetminus Y$.  Since $\Gamma $ is neat, $\bar Y$
is a compact manifold with corners of dimension $N$ that is homotopy
equivalent to $Y$.  Thus
\[
H^{*} (\Gamma ;\Z )=H^{*} (Y;\Z )= H^{*} (\bar Y;\Z )
\] 
and by Lefschetz duality 
\[
H^{*} (\bar Y;\Z )=H_{N-*} (\bar Y, \partial \bar Y; \Z ).
\]

Rather than use $\bar Y$ to compute the cohomology of $\Gamma $, we
will use a certain singular compactification of $Y$, which was
originally constructed by Satake~\cite{satake2, satake1}.  Recall
that $\tilde{C}\subset \bar C$ is the union of $C$ with its rational
boundary components.  We equip $\tilde{C}$ with the \emph{Satake
topology}~\cite{satake2}.  Let $\tilde{X}$ be the quotient of
$\tilde{C}$ by homotheties, and let $\tilde{Y}=\Gamma \backslash
\tilde{X}$.  Then $\tilde{Y}$ is a compact hausdorff space called a
\emph{Satake compactification} of $Y$, and we write $\partial
\tilde{Y} = \tilde{Y}\smallsetminus Y$.

\begin{proposition}\label{comparing.coho}
$H_{*}(\bar Y, \partial \bar Y; \Z )=H_{*}
(\tilde{Y},\partial \tilde{Y};\Z )$. 
\end{proposition}

\begin{proof}
According to Zucker~\cite{zucker}, there is a quotient map $q:\bar
Y\rightarrow \tilde{Y}$ that is the identity on $Y$.  Let $U\subset
\bar Y$ be any collared neighborhood of $\partial \bar Y$; this admits
a deformation retraction to $\partial \bar Y$.  Let $V$ be $q (U)$,
which is a neighborhood of $\partial \tilde{Y}$.  The composition of
$q$ with the retraction yields a retraction of $V$ onto $\partial
\tilde{Y}$.

The result then follows from the sequence of isomorphisms
\[
H_{*} (\bar Y,\partial \bar Y)\xrightarrow{r}H_{*} (\bar Y,U)\xrightarrow{e}
H_{*} (Y,U\cap Y)\xrightarrow{\text{id}}H_{*} (Y,V\cap Y)
\xrightarrow{e}
H_{*} (\tilde{Y},V)\xrightarrow{r}H_{*} (\tilde{Y},\partial \tilde{Y}),
\]
where the maps $r$ are induced by the appropriate retractions, and the
maps $e$ are induced by excision.
\end{proof}

\begin{proposition}\label{using.triangulation}
Under quotienting by homotheties $\pi \colon\tilde{C}\rightarrow \tilde{Y}$,
the \Vor \ cones descend to a finite triangulation of $\tilde{Y}$.
\end{proposition}

\begin{proof}
The \Vor \ cones descend modulo homotheties to give a decomposition of
$\tilde{X}$ into closed simplices; this has the same formal properties
as in \S\ref{fan.stuff}. By Theorem 1 of \cite[p. 113]{amrt}, the
restriction of the Satake topology to each closed simplex
$\tilde{\sigma }$ coincides with the ordinary topology on the simplex.
By a neatness argument like that in Proposition \ref{gam-eq-for-Pi},
$\tilde{\sigma }\cap \gamma \tilde{\sigma } \not = \varnothing $ for
$\gamma \in \Gamma $ implies $\gamma $ fixes $\tilde{\sigma }$
pointwise.  This implies the map $\tilde{X}\rightarrow \tilde{Y}$
restricts to a bijection---hence an embedding---on $\tilde{\sigma }$.
Thus $\tilde{Y}$ has a decomposition into closed simplices with the
same formal properties as in \S\ref{fan.stuff}---that is, a
triangulation of $\tilde{Y}$.  By (3) of \S\ref{gam.admis.section},
the triangulation is finite.
\end{proof}

\subsection{}\label{realization}
Let $\Delta _{k}$ be the standard $k$-simplex, that is 
\[
\Delta _{k} := \Bigl\{(r_{0},\dots, r_{k})\subset \R ^{k+1} \Bigm |
\sum {r_{i}}=1, \quad \hbox{and}\quad r_{i}\geq 0 \Bigr\}.
\]
Let $e_{i}$ be the $i$th vertex of $\Delta _{k}$.  Let $C_{k}
(\tilde{Y})$ be the group of integral singular $k$-chains.  In other
words, $C_{k} (\tilde{Y})$ is the free abelian group generated by all
continuous maps $s\colon \Delta _{k}\rightarrow \tilde{Y}$.

Let $\cones _{k}^{R}$ be the free abelian group generated by the set
of all pointed rational polyhedral cones $\sigma (\xx )$, as $\xx $
varies over all $(k+1)$-tuples of nonzero points in $\tilde{C} (\Q )$.
In general, we shall use boldface to denote complexes built from
polyhedral cones, before taking the quotient by homotheties.  Notice
the shift in degree: $\cones ^{R}_{k}$ is built from
$(k+1)$-dimensional cones.  Cycles in $\cones ^{R}_{k}$ will push
forward to the ``correct'' degree $k$ in $H_{k} (\tilde{Y}, \partial
\tilde{Y})$ and related groups.

The obvious boundary map $\partial \colon\cones _{k}^{R} \rightarrow
\cones _{k-1}^{R} $ makes $\cones _{*}^{R}$ into a chain complex.  By
abuse of notation, we let $\sigma (\xx )$ denote both the pointed
rational polyhedral cone determined by $\xx $ as well as the class of
this cone in $\cones _{k}^{R}$.

Given a chain $\xi = \sum n (\xx )\sigma (\xx ) \in \cones _{k}^{R}$,
where $n (\xx )\in \Z $, we define the \emph{support} of $\xi $ by
\[
\support (\xi ) := \left\{ \sigma (\xx ) \bigm | n (\xx )\not =0 \right\}.
\]

Any chain $\xi$ determines a singular chain $[\xi ]\in C_{k}
(\tilde{Y})$ via the homotheties and mod $\Gamma $.  If $\xx =
(x_{0},\dots ,x_{k})$, the cone $\sigma (\xx )$ induces $s\colon
\Delta _{k}\rightarrow \tilde{Y}$ by $s (e_{i})=\pi (x_{i})$, linear
extension, and taking the quotient modulo $\Gamma $.  The map
$[\phantom{\sigma }]\colon \cones _{*}^{R}\rightarrow C_{*}
(\tilde{Y})$ is a morphism of complexes.

\begin{proposition}\label{liftprop}
Any class in $H_{*} (\tilde{Y},\partial \tilde{Y})$ can be represented by
a chain in $\cones _{*}^{R}$.
\end{proposition}

\begin{proof}
This follows from Proposition~\ref{using.triangulation}.  Any homology
class can be written as the image of a $\Z $-linear combination of
\Vor \ cones, and such cones have a rational pointed structure.
\end{proof}

\begin{definition}\label{real.map}
Given $u\in H_{*} (\tilde{Y},\partial \tilde{Y})$, a chain $\xi \in
\cones _{*}^{R}$ with $[\xi ]=u $ is called a \emph{lift} of~$u$.
\end{definition}

\subsection{}\label{cyclessection}
Let $\zyklen^{R}_{*}$ be the \emph{relative} cycle group
\[
\zyklen ^{R}_{*} := \left\{\xi \in \cones ^{R}_{*} \Bigm | \hbox{$[\partial
\xi ]$ is supported on $\partial \tilde{Y}$.}\right\}
\]
In general, elements of $\zyklen^{R}_{*}$ will not be cycles
with respect to the boundary map in $\cones ^{R}_{*}$: they
will have boundaries mapping to $\partial \tilde{Y}$, and will only be
cycles modulo the $\Gamma $-action. 

Recall that $n$ is the $\Q $-rank of $C$.
\begin{proposition}\label{vanishing}
Let $\xi \in \cones _{k}^{R}$ be supported on cones spanned by cusps.  
\begin{enumerate}
\item If $k=n-1$, then $\xi \in \zyklen  _{k}^{R}$.
\item If $k < n-1$, then the class of $\xi $ is zero in $H_{*}
(\tilde{Y},\partial \tilde{Y})$. 
\end{enumerate}
\end{proposition}

\begin{proof} 
Recall that $C (i)$ is the union of all rational boundary components
of $C$ of $\Q $-rank $\leq i$.  Let $D\subset C (i)$ and $E\subset C
(j)$ be rational boundary components.  According to \cite[Lemma
4]{ash}, we have $D+E\subset C (i+j)$.  Hence if $\xx $ is a $k$-tuple
with $x\in \Xi (C)$ for all $x\in \xx $, then $\sigma (\xx )\subset C
(k)$.

Both statements in the Proposition follow from this containment.  In
(1), $\support (\partial \xi )$ consists of cones that lie in $C
(n-1)$, which under homotheties and $\Gamma $ map to a singular chain
supported on $\partial \tilde{Y}$.  This implies that $\xi $ is a
relative cycle.  In (2), $\support (\xi) $ is a subset of $C (n-1)$,
and itself maps to a singular chain supported on $\partial \tilde{Y}$.
\end{proof}

\begin{example}
Suppose that $G$ is the $\Q $-group with $G (\Q )=SL_{n} (\Q )$, and
let $\xx $ be a rational $n$-tuple with $x\in \Xi (C)$ for all $x\in
\xx $.  Then the class 
\[
[\sigma (\xx )]\in H_{n-1} (\tilde{Y},\partial \tilde{Y}) = H^{N-n+1} (\Gamma )
\]
is a \emph{minimal modular symbol} as in Ash-Rudolph
\cite{ash.rudolph}.
\end{example}

\begin{remark}\label{bs}
Since cycles supported on \Vor \ cones map surjectively onto $H_{*}
(\tilde{Y},\partial \tilde{Y})$, we obtain that
$H_{k}(\tilde{Y},\partial \tilde{Y})$ vanishes for $k<n-1$.
Equivalently, $H^{k} (\Gamma ;\Z )=0$ if $k>N-n+1$.  This is a special
case of a much more general result from~\cite{borel.serre}: the
cohomological dimension of $\Gamma $ is $N-n+1$.
\end{remark}

\subsection{}\label{action.on.coho}
Now we describe the action of ${\HHH}_{\Gamma}$ on $H^{N-*}
(\Gamma )=H_{*} (\tilde{Y},\partial \tilde{Y})$ in the setting of
\S\S\ref{realization}--\ref{cyclessection}.  Choose a Hecke operator
$T_{g}$.  Then we may decompose the double coset $\Gamma g\Gamma $ as
\begin{equation}\label{decomp}
\Gamma  g\Gamma  = \coprod_{s\in S} \Gamma s
\end{equation}
for some set $S\subset G (\Q )$, which is finite since $g\in
\Comm(\Gamma)$.  The Hecke correspondence carries the point
$\Gamma x\in \Gamma \backslash X$ to the finite set of points
$\{\Gamma sx \}_{s\in S}$ in $\Gamma \backslash X$.

Given a class $u\in H_{*} (\tilde{Y},\partial \tilde{Y})$ with lift
$\xi = \sum_{\xx \in A} n (\xx )\sigma (\xx )$, the Hecke operator acts by
\begin{equation}\label{hecke.map}
u\longmapsto \Biggl[ \sum_{\substack{\xx\in A\\ s\in S}} n (\xx )\sigma (s \cdot \xx )\Biggr].
\end{equation}
One can easily show that the image of the map in \eqref{hecke.map} is
a well-defined homology class.

%
%
\section{First Algorithm}\label{first}

As before, let $N$ be the dimension of the symmetric space $X$, and
let $n$ be the $\Q $-rank of $C$.  We want to compute the action of
the Hecke operators on $H^{*} (\Gamma )$, which we identify with the
relative homology $H_{N-*} (\tilde{Y},\partial \tilde{Y})$.
Concretely, we must identify a finite basis\footnote{If $H_{*}
(\tilde{Y},\partial \tilde{Y})$ has torsion, then by ``basis'' we mean
a minimal set of elements generating $H_{*} (\tilde{Y},\partial
\tilde{Y})$ as an abelian group.} of $H_{*} (\tilde{Y},\partial
\tilde{Y})$, and compute the transformation matrix of a given Hecke
operator in terms of this basis.

Let $\zyklen ^{V}_{*}\subset \zyklen ^{R}_{*}$ be the polyhedral
cycles supported on \Vor \ cones.  Let $\zyklen ^{\HHH}_{*}\subset
\zyklen ^{R}_{*}$ be the subgroup of \emph{Hecke images}; that is,
$\xi \in \zyklen ^{\HHH}_{*}$ if and only if $\xi = T (\eta) $ for
some $T\in \HHH_\Gamma $ and $\eta \in \zyklen ^{V}_{*}$.  By
Proposition~\ref{using.triangulation}, $\zyklen ^{V}_{*}$ contains a
finite set of cycles whose homology classes form a basis of $H_{*}
(\tilde{Y},\partial \tilde{Y})$.  Hence to compute the action of a
Hecke operator, it suffices to construct an algorithm that transforms
a cycle in $\zyklen^{\HHH}_{*}$ to cycle in $\zyklen^{V}_{*}$ that
generates the same homology class in $H_{*} (\tilde{Y},\partial
\tilde{Y})$.  In Theorems~\ref{thm1} and~\ref{thm2} we describe an
algorithm that accomplishes this.

\subsection{}\label{notation}
We begin by establishing some constructions appearing in the
algorithm.  Recall that $R (\sigma )$ denotes the set of spanning rays
for a cone $\sigma $.

Let $\sigma \subset \tilde{C}$ be a rational polyhedral cone
satisfying $R (\sigma )\subset \Xi (C)$.  Then $\sigma \cap \V$
denotes the fan obtained by intersecting $\sigma $ with the \Vor \ fan
$\V$.  That is, writing $\{\sigma _{\alpha } \}_{\alpha \in A}$ for
$\sigma \cap \V$, then $\bigcup_{\alpha \in A} \sigma _{\alpha } =
\sigma $, and each $\sigma _{\alpha }$ is the intersection of a \Vor \
cone and some (not necessarily proper) face of $\sigma $.  We call
$\sigma \cap \V$ the \emph{canonical fan} associated to $\sigma $.

Given a cone $\sigma _{\alpha }\in \sigma\cap \V$, let $V_{\alpha}$
denote the \Vor \ cone inducing $\sigma _{\alpha }$.  That is,
$V_{\alpha }$ is the smallest \Vor \ cone containing $\sigma _{\alpha
}$.

\begin{proposition}\label{finite.index.set}
If $\sigma $ is a rational polyhedral cone, then $\sigma \cap \V$ is a
finite fan.
\end{proposition}

\begin{proof}
According to the proof of the main theorem in~\cite{ash}, the
intersection of any rational polyhedral cone with the \Vor \
polyhedron $\Pi $ is cut out by the faces of $\sigma $ and finitely
many supporting hyperplanes of $\Pi $.  Thus $\sigma $ can meet at
most a finite number of \Vor \ cones, and the result follows.
\end{proof}

Let $\xx $ be a rational tuple.  Although $\sigma (\xx )$ and the
cones in $\V$ are simplicial, the canonical fan $\sigma (\xx )\cap \V$
need not be.  Let $F (\xx )$ denote a simplicial refinement of the
canonical fan that does not add any new $1$-cones.

\subsection{}\label{minimal}
We are now ready to present our first algorithm.  In
Theorem~\ref{thm1}, we describe the algorithm on cycles in $\zyklen
_{n-1}^{\HHH}$.  Theorem~\ref{thm2} describes the algorithm on cycles
in $\zyklen _{k}^{\HHH}$ for $k>n-1$.

Let $\xi \in \zyklen _{n-1}^{\HHH}$ be a Hecke image.  Write $\xi $ as
$\sum n (\xx ) \sigma (\xx )$, where $n (\xx ) \in \Z $.  Since
$k=n-1$, Proposition~\ref{vanishing} implies that each $\sigma (\xx )$
is already a cycle.  Hence we may assume without loss of generality
that $\xi =\sigma (\xx )$.

\begin{theorem}\label{thm1}
Let $\xi = \sigma (\xx )\in \zyklen _{n-1}^{\HHH }$.  The following
algorithm constructs a cycle $\xi^{V} \in \zyklen _{n-1}^{V}$ such
that $[\xi ] = [\xi ^{V}]$.
\begin{enumerate}
\item   Construct the canonical fan $\sigma (\xx )\cap \V$,
and construct a simplicial refinement $F (\xx )$ that does not add any
new $1$-cones.  
\item   For each $1$-cone $\rho _{\beta }\in F (\xx )$,
select a cusp $v_{\beta }\in R (V_{\beta })\cap \Xi (C')$.  Here
$V_{\beta }$ is the \Vor \ cone inducing $\rho _{\beta }$, and $C'$ is
the smallest (not necessarily proper) rational boundary component
containing $\rho _{\beta }$.  Let $y_{\beta }$ denote the vertex of $\Pi $
generating $v_{\beta }$. 
\item   For each $n$-cone $\sigma _{\alpha }\in F (\xx )$,
construct a pointed simplicial $n$-cone $\tau (\yy _{\alpha })$ by
taking the $n$ points $\left\{y_{\beta }\mid \rho _{\beta }\subset
\sigma _{\alpha } \right\}$ and ordering them so that the induced
orientation on $\tau (\yy _{\alpha })$ matches the orientation
$\sigma_{\alpha } $ inherits from $\sigma $.
\end{enumerate}
Then $\tau (\yy _{\alpha })\in \zyklen _{n-1}^{V}$, and the desired
cycle is
\[
\xi^{V} = \sum_{\sigma _{\alpha }\in F (\xx )_{n}}\tau (\uu _{\alpha }).
\]
\end{theorem}

\subsection{}\label{correctness.one}
We begin the proof by showing that the choices in Step~2 of Theorem~\ref{thm1}
are possible. 

\begin{lemma}\label{possible}
Let $\rho_{\beta } \in F (\xx )$ be a $1$-cone.  Then the set $R (V_{\beta
})\cap \Xi (C')$ is nonempty.
\end{lemma}

\begin{proof}
By Proposition~\ref{bdy}, the intersection $\Pi ' := C'\cap \Pi $ is a
face of $\Pi $.  The \Vor \ cone $V_{\beta }$ is the cone over a face
$A_{\beta }$ of $\Pi $.  So $V_{\beta}\cap C'$ is the cone over the
face $A_{\beta }\cap \Pi '$ of $\Pi $, which is nonempty since
$V_{\beta }\cap C'\not =\varnothing$.  The set $R (V_{\beta })\cap \Xi
(C')$ consists of the cusps generated by the vertices of $A_{\beta
}\cap \Pi '$.
\end{proof}

Now we construct a cycle that serves as a step between $\sigma (\xx )$
and $\xi^{V} $.

\begin{lemma}\label{xiprime}
The fan $F (\xx )$ induces a cycle $\xi '\in \zyklen _{n-1}^{R}$ with
$[\sigma (\xx )] =[\xi ']$ in $H_{n-1} (\tilde{Y},\partial \tilde{Y})$.
\end{lemma}

\begin{proof}
Essentially, we use $F (\xx )$ to refine $\sigma (\xx )$.
From each $1$-cone $\rho \in F (\xx )$, choose a nonzero $x_{\rho }\in
V (\Q )$.  Then for each $n$-cone $\sigma _{\alpha }\in F (\xx
)_{n}\smallsetminus F (\xx )_{n-1}$, use the points $\left\{x_{\rho }\mid
\rho \subset \sigma _{\alpha }\right\}$ to construct $n$-tuples $\xx
_{\alpha }$ as in Step 3 of Theorem~\ref{thm1}.  These tuples
generate pointed simplicial cones $\sigma (\xx_{\alpha
})$, and we then take 
\[
\xi ' = \sum _{\sigma _{\alpha }\in F (\xx)_{n}} \sigma
(\xx _{\alpha }).
\]
\end{proof}

For the cycle $\xi^{V} $ to lie in $\zyklen _{n-1}^{V}$, we must show that
the cones constructed in Step~3 are \Vor \ cones.

\begin{lemma}\label{vorcone}
Each $\tau (\yy _{\alpha})$ is a \Vor \ cone.
\end{lemma}

\begin{proof}
Let $\sigma _{\alpha }\in F (\xx )_{n}$ and $\tau (\yy _{\alpha })$ be
as in Step 3.  Each cusp $v_{\beta }$ used to construct $\tau (\yy
_{\alpha })$ is the spanning ray of a \Vor \ cone $V_{\beta }$ that
induced a $1$-cone $\rho _{\beta }\subset \sigma _{\alpha }$.  Since
$\rho _{\beta } $ ranges over the $1$-cones with $\rho _{\beta
}\subset \sigma _{\alpha }$, it follows that $V_{\beta}\subseteq
V_{\alpha }$.  Thus all the $v_{\beta }$ lie in $R (V_{\alpha })$.
Since $V_{\alpha}$ is a simplicial cone, the result follows.
\end{proof}

\begin{lemma}\label{eta.lemma}
The class in $H_{n-1} (\tilde{Y},\partial \tilde{Y})$ generated by 
\[
\xi^{V} = \sum_{\sigma _{\alpha }\in F (\xx )_{n}}\tau (\uu _{\alpha })
\]
is equal to
$[\sigma (\xx )]$.
\end{lemma}

\begin{proof}
We will construct a chain $\eta \in \cones _{n}^{R}$ so that $\partial
\eta =\sigma (\xx )-\sum \tau (\uu _{\alpha }) + \mu $, where $\support
(\mu )$ consists of cones lying in $\tilde{C}\smallsetminus C$.  

Consider the fan of simplicial cones $F (\xx )$.  By taking the
intersection of a generic affine hyperplane $H$ with $F (\xx )$, we
obtain a convex union of $(n-1)$-simplices in $H$, which we denote by
$P$.  

Now consider the product $P\times \Delta _{1}$.  This can be realized
as the convex union of a set of simplicial prisms in some affine space
$H'\supset H$.  (See Figure~\ref{constructingP}.)  We subdivide
$P\times \Delta _{1}$ into simplices without adding new vertices, and
call the resulting union $P'$.

\begin{figure}[ht]
\begin{center}
\psfrag{F}{$F (\xx )$}
\psfrag{P}{$P$}
\psfrag{P1}{$P'$}
\psfrag{Dp}{$P\times \Delta _{1}$}
\includegraphics{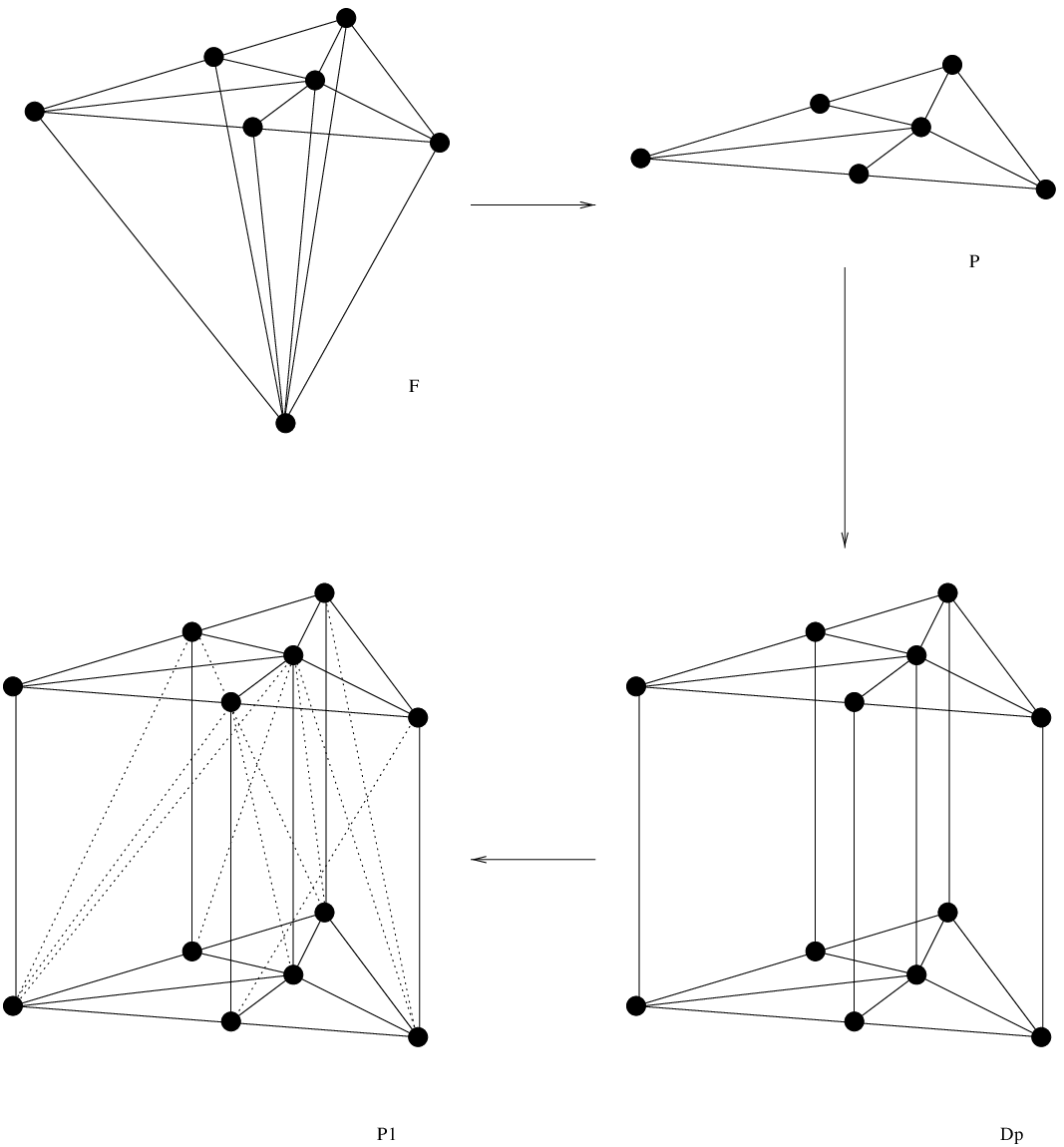}
\end{center}
\caption{\label{constructingP}}
\end{figure}

Next, we linearly map $P'$ to $\tilde{C}$ as follows.  On the upper
face, we take the vertices of $P'$ to the vertices of $\Pi $ that give
the pointed structure to $\sum \tau (\uu _{\alpha })$.  On the lower
face of $P'$, we take the vertices to the points $\{x_{\rho} \}$ in
Lemma~\ref{xiprime} that give the pointed structure to the cycle $\sum
\sigma (\xx _{\alpha })$.  This defines a collection of pointed
$(n+1)$-cones in $\tilde{C}$, and hence defines a chain $\eta $.

From the construction, we have $\partial \eta = \sigma (\xx ) - \sum
\tau (\uu _{\alpha })+ \mu $, where $\mu $ is a chain induced by the
outer sides of $P'$.  To finish the proof, we must show that $\support
(\mu )\subset \tilde{C}\smallsetminus C$.  

To see this, let $\sigma '$ be a maximal proper face of $\sigma (\xx
)$, and let $P (\sigma ')$ be the subset of $P'$ corresponding to
$\sigma '\times \Delta _{1}$.  By Proposition~\ref{vanishing}, we have
$\sigma '\subset \tilde{C}\smallsetminus C$.  In fact, $\sigma '$
determines a unique proper rational boundary component $C'$ such that
$\sigma '\subset C'$.

Any 1-cone in $F (\xx )$ that meets $\sigma '$ is a subset of $\sigma
'$, and hence also lies in $C'$.  This means
that the intersection of $P (\sigma ')$ with the lower face of $P'$
determines a collection of cones lying in $C'$.

Now consider the portion of $P (\sigma ')$ lying in the upper face of
$P'$.  The 1-cones determined by this subset include the spanning rays
of $\sigma '$ as before, as well some cusps used in the construction
of the $\yy _{\alpha }$.  By Proposition~\ref{bdy} and the
restrictions placed on these cusps in Step~2 of Theorem~\ref{thm1},
these cusps also lie in $C'$.  Hence all the cones in the image of $P
(\sigma ')$ lie in $C'$.  Applying this argument to all maximal proper
faces of $\sigma (\xx )$
completes the proof.

\end{proof}

Finally we complete the proof of the theorem.

\begin{proof}[Proof of Theorem~\ref{thm1}] By Lemma~\ref{possible}, we
may choose the cusps in Step~1 without obstruction.
Lemmas~\ref{xiprime} and~\ref{eta.lemma} imply that the cycle $\xi^{V}
$ from Step~3 satisfies $[\xi ]=[\xi ^{V}]$, and Lemma~\ref{vorcone}
implies that $\xi^{V} $ is supported on \Vor \ cones.
\end{proof}
 
\subsection{}\label{nonminimal}
Now we modify our algorithm to work on cycles $\xi \in \zyklen
_{k}^{\HHH}$, where $k>n-1$.  The basic construction is the same, but
there is one crucial difference.  When $k=n-1$, the boundary of any
$\xi = \sum n (\xx )\sigma (\xx )$ automatically lies in $\tilde{C}$.
But for $\xi $ to be a relative cycle when $k>n-1$, its boundary can
meet $\tilde{C}\smallsetminus C$; this part must vanish mod $\Gamma $.
This means that we must be vigilant in our construction to ensure that
the choices in Theorem~\ref{thm1} and
Lemmas~\ref{possible}--\ref{eta.lemma} can be made
equivariantly over $\xi $.

\begin{theorem}\label{thm2}
Let $\xi \in \zyklen _{k}^{\HHH}$ with $k>n-1$.  Then the choices in
Theorem~\ref{thm1} can be made $\Gamma $-equivariantly.  Specifically,
\begin{enumerate}
\item The fans $\{F (\xx )\mid \sigma (\xx )\in \support (\xi )\}$ can
be constructed $\Gamma $-equivariantly, yielding a fan $F (\xi )$.
\item For any two 1-cones $\rho, \rho '\in F (\xi )$ satisfying
$\gamma \cdot \rho= \rho '$ with $\gamma \in \Gamma $, the
corresponding cusps $v$ and $v'$ can be chosen in Step 2 to satisfy $\gamma
\cdot v = v'$.
\end{enumerate}
Furthermore, the construction of $P'$ from $P\times \Delta _{1}$ in
Lemma~\ref{eta.lemma} can be performed $\Gamma $-equivariantly over
all of $\support (\xi )$.  Thus, the algorithm in Theorem~\ref{thm1},
with the modifications above, can be used to construct a cycle $\xi
^{V}\in \zyklen _{k}^{V}$ satisfying $[\xi ]=[\xi^{V}] $.
\end{theorem}

\begin{proof}
The first assertion is the key, for if the fans $F (\xx )$ can be
constructed $\Gamma $-equivariantly, then any construction using them
can be done $\Gamma $-equivariantly by choosing representatives for
the $\Gamma $-orbits.  

Now the intersection of $\support (\xi )$ with $\V$ consists of a
finite set of simplicial cones, so a fortiori there are only a finite
number of such cones modulo $\Gamma $.  This means that we may adapt
the proof of Proposition~\ref{gam-eq-for-Pi} and construct the
refinements using induction on the dimension.
\end{proof}

%
%
\section{Second Algorithm}\label{second}
%
%
In this section we show how to transform cycles in $\zyklen
_{*}^{\HHH}$ to cycles in 
$\zyklen _{*}^{V}$ without explicitly
constructing the fans $F (\xx )$ from \S\ref{notation}.  Instead, we
assume only the existence of an ``oracle'' that answers the following question:
\begin{quotation}
\emph{Given a point $x\in \tilde{C}$, which \Vor \ cone contains $x$?}
\end{quotation}
This oracle for $x\in C$ is the \emph{\Vor \ reduction algorithm}
\cite{voronoi1}\cite{gunn}.  Since the rational boundary components
are self-adjoint homogeneous cones of lower rank, and they receive a
\Vor \ decomposition by $\Pi $ (Proposition~\ref{bdy}), we are
justified in assuming the existence of this oracle on both $C$ and its
boundary components.  Furthermore, in practice the oracle is not
difficult to implement, and requires only slightly more information
than that necessary for cohomology computations with
$\V$.\footnote{See \cite{jaquet} for implementation in the case of
$SL_{2} (\Z )$.}  Accordingly, we define

\begin{definition}\label{sofx}
Given a point $x\in \tilde{C}$, let $S (x)$ be the set of cusps that
span the unique smallest \Vor \ cone containing $x$.  By abuse of
notation, we also write $S (\rho)$ for a 1-cone $\rho $.
\end{definition}

The idea behind our second algorithm is this.  Although it is
expensive to construct the fans $F (\xx )$, it is easy to subdivide a
lift $\xi \in \zyklen _{*}^{\HHH}$ into a cycle $\xi '$ supported on
smaller pointed simplicial cones.  If the top-dimensional cones in
$\xi '$ are small enough, then any two nearby 1-cones $\rho _{1},\rho
_{2}\in \support (\xi ')$ will lie in the same or adjacent cones of $F
(\xx )$.  Thus $S (\rho _{1})\cap S (\rho _{2})$ will be nonempty.
Hence we can hope to choose cusps $v_{i}\in S (\rho _{i})$ and
assemble them into \Vor \ cones as in Theorems~\ref{thm1} and
~\ref{thm2} to build a cycle in $\zyklen _{*}^{V}$.  This is the
approach we take, but our task is complicated by two issues:
\begin{itemize}
\item For efficiency, we want to subdivide as little
as possible.
\item All choices must be made $\Gamma $-equivariantly.
\end{itemize}

To address these issues, we introduce \emph{sufficiently fine
decompositions} of a fan with respect to $\V$
(\S\ref{suff.fine.section}) and the \emph{relative barycentric
subdivision} of a fan with respect to a subfan (\S\ref{bs.section}).
Theorem~\ref{thm3} describes an algorithm that constructs sufficiently
fine decompositions, and Theorem~\ref{thm4} shows how to use such a
decomposition to transform a cycle in $\zyklen _{*}^{\HHH }$ to a
cycle in $\zyklen _{*}^{V}$ in the same homology class.

\subsection{}\label{bs.section}
Let $\Delta = \Delta _{k}$ be the standard $k$-simplex (\S\ref{realization}),
and let $\sets{k}$ be the finite set $\{0,1,\dots ,k \}$.  There is a
bijection between faces of $\Delta $ and subsets of $\sets{k}$:
the subset $I\subset \sets{k}$ corresponds to the convex hull of $\{e_{i}\mid
i\in I\}$.

Recall that the \emph{barycentric subdivision} of $\Delta $ is the
simplicial complex with vertices corresponding to nonempty subsets of
$\sets{k}$, and with $i$-faces corresponding to proper flags of length $i+1$:
\[
\varnothing \subsetneqq I_{0}\subsetneqq \dots \subsetneqq
I_{i}\subsetneqq \sets{k}.
\]
We denote the barycentric subdivision by $\B(\Delta )$.  Using the
concrete description of $\Delta $ given in \S\ref{realization}, we
may realize the vertex of $\B(\Delta )$ corresponding to $I\subset
\sets{k}$ by the point $(\sum _{i\in I} e_{i})/\#I$.

\begin{figure}[ht]
\begin{center}
\includegraphics{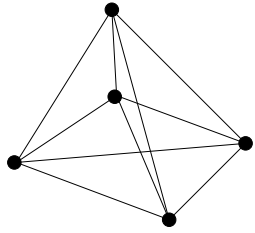}
\end{center}
\caption{\label{starpoint1}The stellar subdivision of $\Delta _{3}$.}
\end{figure}

There is another approach to the barycentric subdivision using the
\emph{stellar subdivision} $\SSS (\Delta )$.  As an abstract simplicial
complex, $\SSS (\Delta )$ is isomorphic to the $(k+1)$-simplex.
In the embedding of $\Delta $ given in \S\ref{realization}, $\SSS
(\Delta )$ is the image of the affine map $\Delta _{k+1}\rightarrow
\Delta _{k}$ that takes
$e_{k+1}\mapsto (\sum _{i\in \sets{k}} e_{i})/ (k+1)$ and is the
identity on the other vertices (Figure~\ref{starpoint1}).  The image of
$e_{k+1}$ is called the \emph{star point}.

\begin{figure}[ht]
\begin{center}
\includegraphics{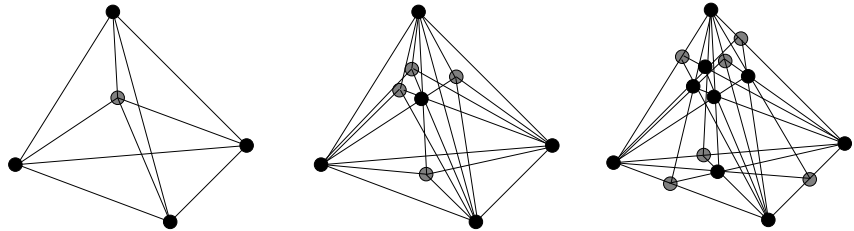}
\end{center}
\caption{\label{starpoint2}Using the stellar subdivision to construct $\B
(\Delta )$.  The grey vertices are added at each step.}
\end{figure}

The barycentric subdivision can now be constructed as follows.  First
construct $\SSS (\Delta )$.  Then stellar subdivide the original
$(k-1)$-faces of $\Delta $, and take the cone of these faces using the
star point of $\SSS (\Delta )$ as the cone point.  Continue
subdividing the lower-dimensional faces of the original $\Delta $,
each time coning with the previously constructed star points.  After
the $1$-faces have been subdivided, the result is $\B (\Delta)$
(Figure~\ref{starpoint2}).

\begin{figure}[ht]
\begin{center}
\includegraphics{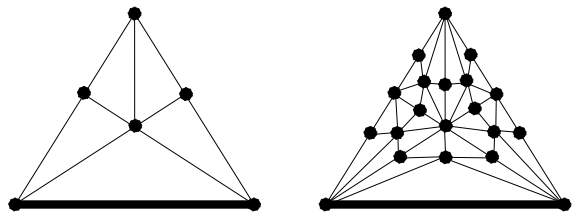}
\end{center}
\caption{\label{relative}$\B ^{1} (\Delta ,\cup
\Delta_{\alpha })$ and $\B ^{2} (\Delta ,\cup
\Delta_{\alpha })$.  The
heavy line represents $\cup \Delta _{\alpha }$.}
\end{figure}

Now let $\cup \Delta_{\alpha }$ be a closed union of proper faces of
$\Delta$.  We define the \emph{relative barycentric subdivision} $\B
(\Delta, \cup \Delta_{\alpha })$ as follows.  We construct the
sequence of stellar subdivisions above, but do not subdivide any
simplices in $\cup \Delta_{\alpha }$ (Figure~\ref{relative}).

Both the barycentric subdivision and the relative barycentric
subdivision can be iterated, and we denote the resulting simplicial
complexes by $\B^{i} (\Delta)$ respectively $\B ^{i} (\Delta ,\cup
\Delta_{\alpha })$. 

Let $\sigma (\xx )$ be a pointed simplicial cone.  Then the
barycentric subdivision of $\sigma (\xx )$ is the fan $\B (\sigma (\xx
))$ constructed by barycentrically subdividing the simplex $\Delta$
spanned by the tuple $\xx $.  We denote the 1-cone on the barycenter
of $\Delta $ by $\beta _{\sigma }$.  We may similarly construct the
relative barycentric subdivision of $\sigma (\xx )$ with respect to a
union of proper faces.  Finally, we may apply either construction to
a fan of pointed simplicial cones by applying it to each simplex
separately.

\begin{remark}\label{dense.remark}
Elementary arguments show that the set of vertices in $\B ^{i} (\Delta
, \cup \Delta _{\alpha })$ becomes dense in $\Delta $ as $i\rightarrow
\infty $, no matter what proper subcomplex $\cup \Delta _{\alpha }$ we
choose.
\end{remark}

\subsection{}\label{suff.fine.section}
We want to formalize how much a given fan must be subdivided for our
purposes.  Recall that $R (\tau )$ denotes the set of spanning rays of
a cone $\tau $, and that if $\rho \subset \tilde{C}$ is a 1-cone, then
$S (\rho )$ denotes the set of spanning rays of the unique smallest
\Vor \ cone containing $\rho $.

\begin{definition}\label{suff.fine}
Let $\sigma \subset \tilde{C}$ be a $d$-dimensional simplicial cone,
and let $\Sigma$ be a finite simplicial fan with $\bigcup_{\tau \in
\Sigma } \tau =\sigma $.  Then the fan $\Sigma $ is said to be a
\emph{sufficiently fine decomposition of $\sigma $} if the following
is true: for every $d$-cone $\tau \in \Sigma $, we have 
\begin{equation}\label{suff.fine.eqn}
S (\tau ):= S(\beta _{\tau })\cap \bigcap _{\rho \subset R(\tau) }S
(\rho )\not = \varnothing.
\end{equation}
More generally, we can speak of a sufficiently fine decomposition of a
set of $d$-cones.  In this case, we mean that each cone in the set has
been refined into a fan whose $d$-cones satisfy \eqref{suff.fine.eqn}.
\end{definition}

\begin{example}\label{2cone}
Suppose that $\sigma $ is a $2$-cone, i.e. $\sigma $ is a cone on an
interval $I$.  Then a partition of $I$ into intervals
$[x_{i},x_{i+1}]$ induces a sufficiently fine decomposition of $\sigma
$ if $x_{i}$ and $x_{i+1}$ lie in adjacent or the same \Vor \ cones.  (The
condition on barycenters is automatic in this case.)  This is
equivalent to the notion of sufficiently fine for modular symbols for
$\Q $-rank one groups presented in~\cite{gunn}.
\end{example}

\subsection{}\label{lemmata}
We want to discuss the relationship between sufficiently fine
decompositions and the canonical fans of \S\ref{notation}.  We begin
by introducing an open covering of a fan.

\begin{definition}\label{open.cover}
Let $F$ be a fan.  Given any $\sigma \in F$, let $U_{\sigma } (F)$ be the
open star of $\sigma $ in $F$.  In other words,
\[
U_{\sigma } (F) = \bigcup_{\substack{\sigma \supseteq \sigma '\\
\sigma '\in F}} \Int \sigma '.
\]
Let $\UUU (F)$ be $\{U_{\sigma } (F)\mid \sigma \in F \}$.
\end{definition}

Notice that if $\sigma \subset \sigma '$, then $U_{\sigma } (F)
\supset U_{\sigma '} (F)$.  Also, if $F'\subset F$ is a subfan (a
subset of $F$ that is also a fan), then 
\begin{equation*}
U_{\sigma } (F') = U_{\sigma } (F) \cap F'.
\end{equation*}
In general the sets in $\UUU (F)$ are not convex, although they are
star-shaped.

Let $F (\sigma )$ be the canonical fan $\sigma \cap \V$, with no
simplicial refinement.

\begin{proposition}\label{big.lemma}
Let $\sigma $ be a simplicial $d$-cone, and let $\Sigma $ be a
simplicial fan that is a refinement of $\sigma $.  Let $\tau \in
\Sigma $, and let $\beta _{\tau }$ be the barycenter of $\tau $.  Then
$\tau $ satisfies \eqref{suff.fine.eqn} if and only if there is a set
$U\in \UUU (F (\sigma ))$ such that $R (\tau )\cup \{\beta _{\tau }
\}\subset U$.
\end{proposition}

\begin{proof}
Assume $\tau $ satisfies \eqref{suff.fine.eqn}.  Then there is a \Vor
\ cone $V_{\alpha }$ such that $S (\tau )\subset R (V_{\alpha })$.
This \Vor \ cone induces $\sigma _{\alpha }\in F (\sigma )$.  The open
star $U_{\sigma _{\alpha }}$ of $\sigma _{\alpha }$ in the canonical
fan is the intersection of $\sigma $ with the open star of $V_{\alpha
}$ in $\V$.  Then $R (\tau )\cup \{\beta _{\tau } \}$ is contained in
$U_{\sigma _{\alpha }}$.  The converse is obtained by simply reversing
this argument.

\end{proof}

A consequence of the preceding Proposition is that sufficiently fine
decompositions exist for any rational simplicial cone in $\tilde{C}$.

\begin{proposition}\label{exist.suff.fine} 
Let $\sigma $ be a rational simplicial $d$-cone.
Then the fan $\B ^{i} (\sigma )$ gives a sufficiently fine
decomposition of $\sigma $ for $i>\!>0$.
\end{proposition}

\begin{proof}
By Proposition~\ref{big.lemma}, any fan subordinate to $\UUU (F
(\sigma ))$ gives a sufficiently fine decomposition of $\sigma $.
Hence the statement follows immediately from the fact that $\UUU
(F (\sigma ))$ is finite.  
\end{proof}

\subsection{}\label{algorithm}
We now describe an algorithm that constructs sufficiently fine
decompositions of a given set of cones.  Recall that if $\Sigma $ is a set of
cones, the we let $\Sigma _{d} = \left\{\sigma \in \Sigma _{d}\mid
\dim \sigma \leq d \right\}$.

\begin{figure}[ht]
\begin{center}
\includegraphics{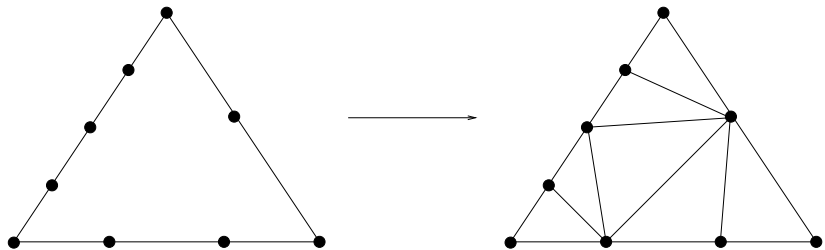}
\end{center}
\caption{An affine slice of the extension of $\Sigma '_{1}$ to $\bar \Sigma $.\label{nonew}}
\end{figure}

\begin{theorem}\label{thm3}
Let $\xi \in \zyklen _{k}^{\HHH}$, and let $\Sigma =\support (\xi )$.
Then the following algorithm constructs a sufficiently fine
decomposition $\Sigma '$ of $\Sigma $:
\begin{enumerate}
\item Set $\Sigma '_{1} = \Sigma _{1}$, and set $j=1$.
\item Assume that $\Sigma' _{j}$ has been constructed for
$1\leq j\leq d$, where $d<k$.  Let $\bar \Sigma $ be any $\Gamma
$-equivariant extension of $\Sigma '_{j}$ to $\Sigma _{j+1}$ without
adding new 1-cones (see Figure~\ref{nonew}).  Then for some $i>0$, the
fan given by the relative barycentric subdivision $\B ^{i}(\bar
\Sigma, \Sigma '_{j})$ is sufficiently fine, and can be taken for
$\Sigma '_{j+1}$.
\item If $j+1=k$, terminate.  Otherwise, increment $j$ and return to Step 2.
\end{enumerate}
\end{theorem}

\begin{proof}
We must verify two facts.  First we claim that $\bar \Sigma $ can be
constructed $\Gamma$-equivariantly.  This follows because each of the
sets $\Sigma _{j+1}\smallsetminus \Sigma _{j}$ is finite modulo
$\Gamma $, and so we may subdivide a set of representatives and
translate by $\Gamma $, as in Proposition~\ref{gam-eq-for-Pi}.

Now we claim that $\B ^{i}(\bar \Sigma, \Sigma '_{j})$ will be
sufficiently fine for $i>\!\!>0$.  If we were performing the usual
barycentric subdivision, then by Proposition~\ref{exist.suff.fine} we
would succeed.  We must show using the relative barycentric
subdivision with respect to $\bar \Sigma_{j} $ poses no obstruction.

So let $\tau \in \Sigma '_{j}$ be a $j$-cone that appears in the
refinement of a given $\sigma \in \Sigma _{j}$.  Let $\beta_\tau $ be
the barycenter of $\tau $.  Since $\Sigma' _{j}$ is sufficiently fine,
Proposition~\ref{big.lemma} implies that there exists a set $U_{\alpha
} \in \UUU (F (\sigma ))$, the open star $\sigma _{\alpha }\subset
F (\sigma )$, such that
\[
R (\tau )\cup \{\beta _{\tau }\}\subset U_{\alpha }.
\]
(See Figure~\ref{side}.
In this and the following figure, we depict a generic affine   
slice of these objects.)

\begin{figure}[ht]
\psfrag{s}{$\sigma $}
\psfrag{t}{$\tau $}
\psfrag{b}{$\beta_\tau $}
\psfrag{Uv}{$U_{\alpha }$}
\begin{center}
\includegraphics{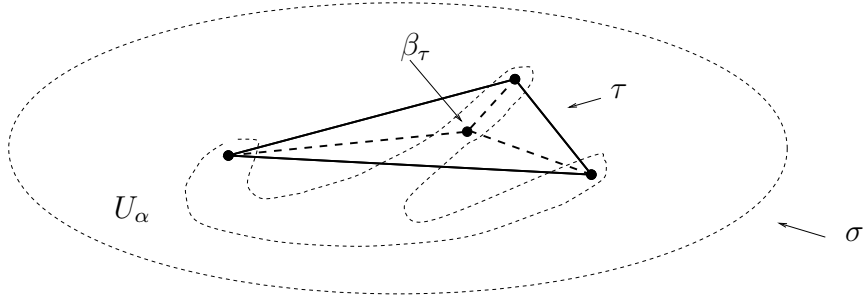}
\end{center}
\caption{$U_{\alpha }$ and $\tau $ in $\sigma $.\label{side}}
\end{figure}

Suppose that $\sigma $ appears as a face of a simplex $\sigma '\in
\Sigma _{j+1}$.  Then the canonical fan $F (\sigma )$ is a subfan of
$F (\sigma ')$, which implies $\sigma _{\alpha }\subset F (\sigma ')$.
Moreover, $U_{\alpha }$ is the intersection of $\sigma $ with the set
$U_{\alpha } (F (\sigma ')) \subset \UUU ( F (\sigma '))$.  We can
picture $U_{\alpha }(F (\sigma '))$ as a ``thickening'' of $U_{\alpha
}$ into $\sigma '$.

Now extend $\Sigma _{j}'$ to $\bar \Sigma $, and begin iterating the
relative barycentric subdivision.  Away from $\partial \sigma '$,
the new $1$-cones will fill $\sigma' $ densely, and by
Proposition~\ref{exist.suff.fine} the decomposition will become
sufficiently fine there.  The only possible problem is near $\partial
\sigma '$.  Over each $\tau \in \partial \sigma '$, the new $1$-cones
of $\B ^{i}(\bar \Sigma, \Sigma '_{j})$ will converge to $\beta _{\tau
}$.  In Figure~\ref{tents} we show the affine slices of these
$1$-cones as grey dots.  By Remark~\ref{dense.remark}, these new
$1$-cones come arbitrarily close to $\beta _{\tau }$ and eventually
enter $U_{\alpha }(F (\sigma '))$.  By Proposition~\ref{big.lemma},
this implies $\B ^{i}(\bar \Sigma, \Sigma '_{j})$ will be sufficiently
fine for some $i>\!>0$.

\begin{figure}[ht]
\psfrag{s}{$\sigma $}
\psfrag{s'}{$\sigma '$}
\psfrag{t}{$\tau $}
\psfrag{Uv}{$U_{\alpha }$}
\begin{center}
\includegraphics{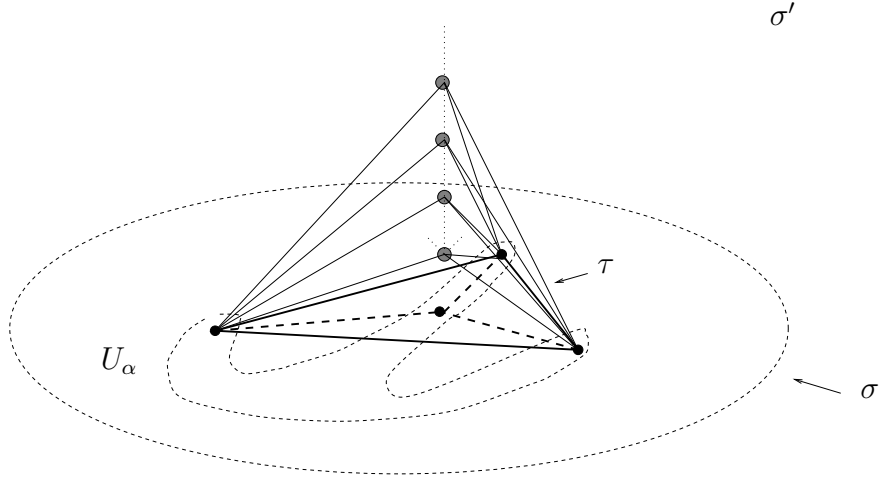}
\end{center}
\caption{The 1-cones converge to $\beta _{\tau }$ and enter $U_{\alpha } (F
(\sigma' ))\supset U_{\alpha }$.\label{tents}}
\end{figure}

\end{proof}

\begin{remark}\label{about.suff.fine}
One might think that the set $S (\beta _{\tau}) $ could be omitted
from \eqref{suff.fine.eqn} in the definition of sufficiently fine
decomposition.  However, since the $U_{\alpha }\in \UUU (F (\sigma ))$
are in general non-convex, the fact that $R (\tau )\subset U_{\alpha
}$ does not imply $\beta _{\tau }\subset U_{\alpha }$.  We need
$\beta_{\tau } \in U_\alpha $ in the proof above to ensure that
sequence of $1$-cones in Figure~\ref{tents} enters $U_{\alpha } (F
(\sigma '))$.
\end{remark}

\subsection{}\label{algorithm.2}
We conclude the paper by showing how a sufficiently fine decomposition
as in Theorem~\ref{thm3} may be used to transform a cycle in $\zyklen
_{*}^{\HHH }$ to a cycle in $\zyklen _{*}^{V}$ giving the same homology
class in $H_{*} (\tilde{Y}, \partial\tilde{Y})$.

\begin{figure}[ht]
\psfrag{sx}{$\tau $}
\psfrag{r1}{$\rho _{1}$}
\psfrag{r3}{$\rho _{2}$}
\psfrag{r2}{$\rho _{3}$}
\begin{center}
\includegraphics{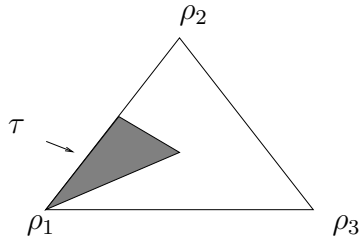}
\end{center}
\caption{$\tau $ corresponds to the
ordered triple $(\rho _{1},\rho _{2},\rho _{3})$.\label{assembling}}
\end{figure}

\begin{theorem}\label{thm4}
Let $\xi =\sum n (\xx )\sigma (\xx )\in \zyklen _{k}^{\HHH }$ be a
Hecke image, and let $F (\xi )$ be a $\Gamma$-equivariant sufficiently
fine decomposition of $\support (\xi )$.  Let $\B (F (\xi ))$
be the barycentric subdivision of $F (\xi )$.  Then the following
algorithm constructs a cycle $\xi ^{V}\in \zyklen _{k}^{V}$ satisfying
$[\xi ] = [\xi ^{V}]$.
\begin{enumerate}
\item Each 1-cone $\rho_{I} \in \support (\B (\xi))$ is the barycenter
of a set $\{\rho _{i} \}_{i \in I}$ of $1$-cones of $F (\xi )$.  For
each $\rho _{I}$, choose a cusp $v_{I}\in \bigcap _{i\in I}S (\rho
_{i})$.  Make these choices $\Gamma $-equivariantly over $\B (F(\xi
))$.  Let $y_{I}$ be the vertex of $\Pi $ generating $v_{I}$.
\item Each $k$-cone $\tau _{\alpha }\in \support (\B (\xi))$
corresponds to a flag
\[
\{\rho _{1} \}\subset \{\rho _{1},\rho _{2} \}\subset \dots \subset
\{\rho _{1},\dots ,\rho _{k} \} 
\] 
where each $\rho _{i}$ is a $1$-cone in $\support (F (\xi ))$
(cf. Figure~\ref{assembling}).  For
each such flag, assemble the points 
\[
y_{\{1 \}}, y_{\{1,2 \}}, \dots , y_{\{1,\dots, k \}}
\]
into a $k$-tuple $\yy _{\alpha }$, using the orientation on $\tau
_{\alpha }$.  Form rational pointed cones $\tau (\yy _{\alpha })$
using these tuples.
\end{enumerate}
Then the desired cycle is 
\[
\xi ^{V}=\sum _{\tau _\alpha \in \B (F (\xi) )_{k}} \tau (\yy _{\alpha }).
\] 
\end{theorem}

\begin{figure}[ht]
\begin{center}
\psfrag{r1}{$\rho _{1}$}
\psfrag{r2}{$\rho _{2}$}
\psfrag{r3}{$\rho _{3}$}
\psfrag{v1}{$v_{\{1 \}}$}
\psfrag{v12}{$v_{\{1,2 \}}$}
\psfrag{v123}{$v_{\{1,2,3 \}}$}
\psfrag{b}{$\beta $}
\includegraphics{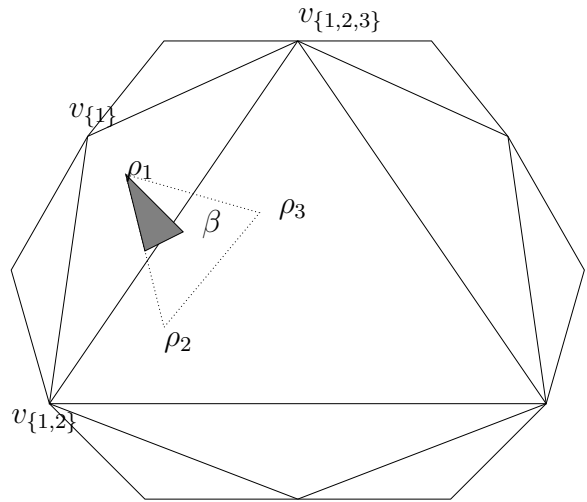}
\caption{An affine slice of part of $\V$ in the case $\Gamma \subset
SL_{2} (\Z )$\label{sl2}.  The shaded triangle is $\tau (\yy _{\alpha })$.}
\end{center}
\end{figure}

\begin{proof}
Let $I$ be a set indexing a $1$-cone $\rho _{I}\in \support (\B (F
(\xi )))$.  Then the definition of sufficiently fine implies that the
set $\bigcap _{i\in I}S (\rho _{i})$ is nonempty.  Hence we may select
a cusp $v_{I}$.  Also, the cusps $v_{I}$ can be chosen $\Gamma
$-equivariantly over all of $\support (\B (F (\xi )))$ by applying a
neatness argument similar to Proposition~\ref{gam-eq-for-Pi}.

Next, that $\xi ^{V}$ is homologous to $\xi $ follows from arguments
identical to those presented in Lemmas~\ref{xiprime}
and~\ref{eta.lemma}.  

What remains to be shown is that the cones $\tau (\yy _{\alpha })$ are
\Vor \ cones.  This follows since the cusps 
\[
v_{\{1 \}}, v_{\{1,2 \}}, \dots , v_{\{1,\dots, k \}}
\]
are all spanning rays of the \Vor \ cone containing $\rho _{1}$.  (See
Figure~\ref{sl2} for an example when $k=2$ and $\Gamma \subset SL_{2}
(\Z )$.)  Since this \Vor \ cone is simplicial, this completes the
proof of the theorem.
\end{proof}

%
%
\bibliographystyle{amsplain}
\bibliography{sahc}

\end{document}